\documentclass{amsproc}
\usepackage{amsfonts}

\newtheorem{theorem}{Theorem}

\newtheorem{definition}[theorem]{Definition}

\newtheorem{lemma}[theorem]{Lemma}

\input{tcilatex}

\begin{document}
\title[ON THE DYNAMICS OF SOLUTIONS OF A RATIONAL DIFFERENCE EQUATION VIA
GENERALIZED TRIBONACCI NUMBERS]{ON THE DYNAMICS OF SOLUTIONS OF A RATIONAL
DIFFERENCE EQUATION VIA GENERALIZED TRIBONACCI NUMBERS}
\thanks{}
\author[\.{I}nci Okumu\c{s}, Y\"{u}ksel Soykan]{\.{I}nci Okumu\c{s}, Y\"{u}%
ksel Soykan}
\maketitle

\begin{center}
\textsl{Zonguldak B\"{u}lent Ecevit University, Department of Mathematics, }

\textsl{Art and Science Faculty, 67100, Zonguldak, Turkey }

\textsl{e-mail: inci\_okumus\_90@hotmail.com \ (corresponding author)}

\textsl{\ yuksel\_soykan@hotmail.com}

\textbf{Abstract}
\end{center}

In this study, we investigate the form of solutions, stability character and
asymptotic behavior of the following rational difference equation%
\begin{equation*}
x_{n+1}=\frac{\gamma }{x_{n}\left( x_{n-1}+\alpha \right) +\beta }\text{, \
\ }n=0,1,...\text{,}
\end{equation*}%
where the inital values $x_{-1}$ and $x_{0}$ and $\alpha $, $\beta $ and $%
\gamma $ with $\gamma \neq 0$ are nonnegative real numbers. Its solutions
are associated with generalized Tribonacci numbers.



\section{Introduction}

Difference equations and systems of difference equations are great
importance in the field of mathematics as well as in other sciences. The
applications of the theory of difference equations appear as discrete
mathematical models of many phenomena such as in biology, economics,
ecology, control theory, physics, engineering, population dynamics and so
forth. Recently, there has been a growing interest in the study of finding
closed-form solutions of difference equations and systems of difference
equations. Some of the forms of solutions of these equations are
representable via well-known integer sequences such as Fibonacci numbers,
Pell numbers, Lucas numbers and Padovan numbers. There are many papers on
such these studies from several authors [\ref{Tollu2013}-\ref{Ocalan2019}].
For example, in [\ref{Bacani2015}], Bacani and Rabago studied the behavior
of solutions of the following nonlinear difference equations%
\begin{equation}
x_{n+1}=\frac{q}{p+x_{n}^{v}}\text{ \ and \ }y_{n+1}=\frac{q}{-p+y_{n}^{v}}%
\text{,}  \label{equ:bsbsbs}
\end{equation}%
where $p$, $q\in 
\mathbb{R}
^{+}$ and $v\in 
\mathbb{N}
$ such that their solutions are associated with Horadam numbers.

Then, El-Dessoky in [\ref{Eldessoky2017}] dealt with the following
difference equation%
\begin{equation}
x_{n+1}=ax_{n}+\frac{\alpha x_{n}x_{n-l}}{\beta x_{n}+\gamma x_{n-k}}\text{,
\ }n=0,1,...\text{,}  \label{equ:kubra}
\end{equation}%
where the parameters $\alpha $, $\beta $, $\gamma $ and $a$ and the initial
conditions $x_{-t}$, $x_{-t+1}$, \ \ ,$x_{-1}$ and $x_{0}$ where $t=\max
\left\{ l,k\right\} $ are positive real numbers. He introduced the explicit
formula of solutions of some special cases of Eq. (\ref{equ:kubra}) via
Fibonacci numbers.

Next, in [\ref{Stevic2018}], Stevic et al. the following nonlinear
second-order difference equation%
\begin{equation}
x_{n+1}=a+\frac{b}{x_{n}}+\frac{c}{x_{n}x_{n-1}}\text{, \ }n\in 
\mathbb{N}
_{0}\text{,}  \label{equ:kjhggtqw}
\end{equation}%
in which parameters $a$, $b$, $c$ and the initial values $x_{-1}$ and $x_{0}$
are complex numbers such that $c\neq 0$. Next, they used the following
change of variables%
\begin{equation*}
x_{n}=\frac{y_{n}}{y_{n-1}}\text{,}
\end{equation*}%
and obtained the following third-order linear difference equation with
constant coefficients%
\begin{equation*}
y_{n+1}=ay_{n}+by_{n-1}+cy_{n-2}\text{.}
\end{equation*}%
After, they introduced the representation formula of every solution of Eq.(%
\ref{equ:kjhggtqw}) with generalized Padovan number.

Alotaibi et al. in [\ref{Alotaibi2018}] considered the following systems of
difference equations%
\begin{equation}
x_{n+1}=\frac{y_{n}y_{n-2}}{x_{n-1}+y_{n-2}}\text{, \ }y_{n+1}=\frac{%
x_{n}x_{n-2}}{\pm y_{n-1}\pm x_{n-2}}\text{, \ }n=0,1,...\text{,}
\label{equ:brtg}
\end{equation}%
where the initial conditions $x_{-2}$, $x_{-1}$, $x_{0}$, $y_{-2}$, $y_{-1}$%
, $y_{0}$ are arbitrary positive real numbers. they analyzed the solutions
of the systems (\ref{equ:brtg}) such that their solutions are associated
with Fibonacci numbers.

In this paper, we study the following difference equation%
\begin{equation}
x_{n+1}=\frac{\gamma }{x_{n}\left( x_{n-1}+\alpha \right) +\beta }\text{, \
\ }n=0,1,...\text{,}  \label{equ:gutreb}
\end{equation}%
where the inital values $x_{-1}$ and $x_{0}$ are arbitrary nonzero real and
the parameters $\alpha $, $\beta $ and $\gamma $ are nonnegative real
numbers with $\gamma \neq 0$.

\section{Preliminaries}

\subsection{Linearized stability}

Let $I$ be some interval of real numbers and let $f:I^{k+1}\rightarrow I$ be
a continuously differentiable function. A difference equation of order $%
(k+1) $ is an equation of the form 
\begin{equation}
x_{n+1}=f(x_{n},x_{n-1},...,x_{n-k}),\text{ \ \ \ }n=0,1,...\text{.}
\label{equ:dsfsgsoyle}
\end{equation}

A solution of Eq.(\ref{equ:dsfsgsoyle}) is a sequence $\{x_{n}\}_{n=-k}^{%
\infty }$ that satisfies Eq.(\ref{equ:dsfsgsoyle}) for all $n\geq -k$.

\begin{definition}
A solution of Eq.(\ref{equ:dsfsgsoyle}) that is constant for all $n\geq -k$
is called an $equilibrium$ $solution$ of Eq.(\ref{equ:dsfsgsoyle}). If 
\begin{equation*}
x_{n}=\overline{x}\text{,\ for all }n\geq -k
\end{equation*}%
is an equilibrium solution of Eq.(\ref{equ:dsfsgsoyle}), then $\overline{x}$
is called an $equilibrium$ $point$, or simply an $equilibrium$ of Eq.(\ref%
{equ:dsfsgsoyle})..
\end{definition}

\begin{definition}[Stability]
Let $\overline{x}$ an equilibrium point of Eq.(\ref{equ:dsfsgsoyle}).

\begin{description}
\item[(a)] An equilibrium point $\overline{x}$ of Eq.(\ref{equ:dsfsgsoyle})
is called locally stable if, for every $\varepsilon >0$; there exists $%
\delta >0$ such that if $\{x_{n}\}_{n=-k}^{\infty }$ is a solution of Eq.(%
\ref{equ:dsfsgsoyle}) with 
\begin{equation*}
\left \vert x_{-k}-\overline{x}\right \vert +\left \vert x_{1-k}-\overline{x}%
\right \vert +...+\left \vert x_{0}-\overline{x}\right \vert <\delta ,
\end{equation*}%
then%
\begin{equation*}
\left \vert x_{n}-\overline{x}\right \vert <\varepsilon ,\text{ \ for all }%
n\geq -k\text{.}
\end{equation*}

\item[(b)] An equilibrium point $\overline{x}$ of Eq.(\ref{equ:dsfsgsoyle})
is called locally asymptotically stable if, it is locally stable, and if in
addition there exists $\gamma >0$ such that if $\{x_{n}\}_{n=-k}^{\infty }$
is a solution of Eq.(\ref{equ:dsfsgsoyle}) with 
\begin{equation*}
\left \vert x_{-k}-\overline{x}\right \vert +\left \vert x_{-k+1}-\overline{x%
}\right \vert +...+\left \vert x_{0}-\overline{x}\right \vert <\gamma ,
\end{equation*}%
then we have 
\begin{equation*}
\lim_{n\rightarrow \infty }x_{n}=\overline{x}.
\end{equation*}

\item[(c)] An equilibrium point $\overline{x}$ of Eq.(\ref{equ:dsfsgsoyle})
is called a global attractor if, for every solution $\{x_{n}\}_{n=-k}^{%
\infty }$ of Eq.(\ref{equ:dsfsgsoyle}), we have 
\begin{equation*}
\lim_{n\rightarrow \infty }x_{n}=\overline{x}.
\end{equation*}

\item[(d)] An equilibrium point $\overline{x}$ of Eq.(\ref{equ:dsfsgsoyle})
is called globally asymptotically stable if it is locally stable, and a
global attractor.

\item[(e)] An equilibrium point $\overline{x}$ of Eq.(\ref{equ:dsfsgsoyle})
is called unstable if it is not locally stable.
\end{description}
\end{definition}

Suppose that the function $f$ is continuously differentiable in some open
neighborhood of an equilibrium point $\overline{x}.$ Let 
\begin{equation*}
q_{i}=\frac{\partial f}{\partial u_{i}}(\overline{x},\overline{x},...,%
\overline{x}),\text{ \ for }i=0,1,...,k
\end{equation*}%
denote the partial derivative of $f(u_{0},u_{1},...,u_{k})$ with respect to $%
u_{i}$ evaluated at the equilibrium point $\overline{x}$ of Eq.(\ref%
{equ:dsfsgsoyle}).

\begin{definition}
The equation 
\begin{equation}
y_{n+1}=q_{0}y_{n}+q_{1}y_{n-1}+...+q_{k}y_{n-k},\text{ }n=0,1,...
\label{equ:sdfgshdgf}
\end{equation}%
is called the linearized equation of Eq.(\ref{equ:dsfsgsoyle}) about the
equilibrium point $\overline{x}$, and the equation 
\begin{equation}
\lambda ^{k+1}-q_{0}\lambda ^{k}-...-q_{k-1}\lambda -q_{k}=0
\label{equ:sdjhgsjdgsdg}
\end{equation}%
is called the characteristic equation of Eq.(\ref{equ:sdfgshdgf}) about $%
\overline{x}$.
\end{definition}

\begin{theorem}[The Linearized Stability Theorem]
\label{theorem:tugbaokumus}Assume that the function $f$ is a continuously
differentiable function defined on some open neighborhood of an equilibrium
point $\overline{x}.$ Then the following statements are true:

\begin{description}
\item[(a)] When all the roots of characteristic equation (\ref%
{equ:sdjhgsjdgsdg}) have absolute value less than one, then the equilibrium
point $\overline{x}$ of Eq.(\ref{equ:dsfsgsoyle}) is locally asymptotically
stable.

\item[(b)] If at least one root of characteristic equation (\ref%
{equ:sdjhgsjdgsdg}) has absolute value greater than one, then the
equilibrium point $\overline{x}$ of Eq.(\ref{equ:dsfsgsoyle}) is unstable.
\end{description}
\end{theorem}

Moreover, the equilibrium point $\overline{x}$ of Eq.(\ref{equ:dsfsgsoyle})
is called $hyperbolic$ if no root of characteristic equation (\ref%
{equ:sdjhgsjdgsdg}) has absolute value equal to one. If there exists a root
of characteristic equation (\ref{equ:sdjhgsjdgsdg}) with absolute value
equal to one, then the equilibrium $\overline{x}$ is called $nonhyperbolic$.

An equilibrium point $\overline{x}$ of Eq.(\ref{equ:dsfsgsoyle}) is called a 
$repeller$ if all roots of characteristic equation\ (\ref{equ:sdjhgsjdgsdg})
have absolute value greater than one.

An equilibrium point $\overline{x}$ of Eq.(\ref{equ:dsfsgsoyle}) is called a 
$saddle$ if one of the roots of characteristic equation (\ref%
{equ:sdjhgsjdgsdg}) is greater and another is less than one in absolute
value.

\begin{theorem}[Clark Theorem]
\label{theorem:lkjhgfvcnm} Assume that $q_{0},q_{1},...,q_{k}$ are real
numbers such that 
\begin{equation*}
\left\vert q_{0}\right\vert +\left\vert q_{1}\right\vert +...+\left\vert
q_{k}\right\vert <1
\end{equation*}%
Then all roots of Eq.(\ref{equ:sdjhgsjdgsdg}) lie inside the unit disk.
\end{theorem}

\subsection{Generalized Tribonacci numbers}

First, from [\ref{Soykan2019a}], consider the generalized Tribonacci
sequence $\left\{ V_{n}\right\} _{n=0}^{\infty }$ defined by the recurrent
relation%
\begin{equation}
V_{n+3}=rV_{n+2}+sV_{n+1}+tV_{n}\text{, \ }n\in 
\mathbb{N}
\text{,}  \label{equ:rewqsaczxd}
\end{equation}%
where the constant coefficients $r$, $s$, $t$ are real numbers and the
special initial conditions%
\begin{equation*}
V_{0}=0\text{, }V_{1}=1\text{, }V_{2}=r\text{.}
\end{equation*}%
The sequence $\{V_{n}\}_{n=0}^{\infty }$ can be extended to negative
subscripts by defining%
\begin{equation*}
V_{-n}=-\frac{s}{t}V_{-(n-1)}-\frac{r}{t}V_{-(n-2)}+\frac{1}{t}V_{-(n-3)}
\end{equation*}%
for $n=1,2,3,...$ when $t\neq 0.$ Hereby, recurrence (\ref{equ:rewqsaczxd})
holds for all integer $n$.

As $\{V_{n}\}_{n=0}^{\infty }$ is a third order recurrence sequence
(difference equation), it's characteristic equation is%
\begin{equation}
x^{3}-rx^{2}-sx-t=0\text{,}  \label{chaequ:jhndf}
\end{equation}%
whose roots are%
\begin{eqnarray*}
\varphi &=&\varphi (r,s,t)=\frac{r}{3}+A+B \\
\chi &=&\chi (r,s,t)=\frac{r}{3}+\omega A+\omega ^{2}B \\
\psi &=&\psi (r,s,t)=\frac{r}{3}+\omega ^{2}A+\omega B
\end{eqnarray*}%
where%
\begin{eqnarray*}
A &=&\left( \frac{r^{3}}{27}+\frac{rs}{6}+\frac{t}{2}+\sqrt{\Delta }\right)
^{1/3},\text{ }B=\left( \frac{r^{3}}{27}+\frac{rs}{6}+\frac{t}{2}-\sqrt{%
\Delta }\right) ^{1/3} \\
\Delta &=&\Delta (r,s,t)=\frac{r^{3}t}{27}-\frac{r^{2}s^{2}}{108}+\frac{rst}{%
6}-\frac{s^{3}}{27}+\frac{t^{2}}{4}, \\
\omega &=&\frac{-1+i\sqrt{3}}{2}=\exp (2\pi i/3)
\end{eqnarray*}%
Notice that we get the following identities%
\begin{eqnarray*}
\varphi +\chi +\psi &=&r, \\
\varphi \chi +\varphi \psi +\chi \psi &=&-s, \\
\varphi \chi \psi &=&t.
\end{eqnarray*}%
From now on, we assume that $\Delta (r,s,t)>0,$ so that the Eq.(\ref%
{equ:rewqsaczxd}) has one real ($\varphi $) and two non-real solutions with
the latter being conjugate complex. Therefore, in this case, it is widely
known that generalized Tribonacci numbers can be stated, for all integers $%
n, $ using Binet's formula%
\begin{equation}
V_{n}=\frac{\varphi ^{n+1}}{(\varphi -\chi )(\varphi -\psi )}+\frac{\chi
^{n+1}}{(\chi -\varphi )(\chi -\psi )}+\frac{\psi ^{n+1}}{(\psi -\varphi
)(\psi -\chi )}\text{.}  \label{equat:mnopcvbedcxzsa}
\end{equation}

We can present Binet's formula of the generalized Tribonacci numbers for the
negative subscripts:\ for $n=1,2,3,...$ we have%
\begin{eqnarray*}
V_{-n} &=&\frac{\varphi ^{2}-r\varphi -s}{t}\frac{\varphi ^{2-n}}{(\varphi
-\chi )(\varphi -\psi )}+\frac{\chi ^{2}-r\chi -s}{t}\frac{\chi ^{2-n}}{%
(\chi -\varphi )(\chi -\psi )} \\
&&+\frac{\psi ^{2}-r\psi -s}{t}\frac{\psi ^{2-n}}{(\psi -\varphi )(\psi
-\chi )}\text{.}
\end{eqnarray*}

\begin{lemma}
Let $\varphi $, $\chi $ and $\psi $ be the roots of Eq.(\ref{chaequ:jhndf}),
suppose that $\varphi $ is a real root with $\max \left \{ \left \vert
\varphi \right \vert ;\left \vert \chi \right \vert ;\left \vert \psi
\right
\vert \right \} =\left \vert \varphi \right \vert $. Then,%
\begin{equation}
\lim_{n\rightarrow \infty }\frac{V_{n+1}}{V_{n}}=\varphi \text{.}
\label{lim:ygbhf}
\end{equation}
\end{lemma}

\textbf{Proof.} Note that there are three cases of the roots, that is when
the roots are all real and distinct, all roots are equal or two roots are
equal, complex conjugate. We will only proof the first case. The proof of
the other two cases of the roots is similar to first one, so it will be
omitted.

If $\varphi $, $\chi $ and $\psi $ are real and distinct then, from Binet's
formula%
\begin{eqnarray*}
\lim_{n\rightarrow \infty }\frac{V_{n+1}}{V_{n}} &=&\lim_{n\rightarrow
\infty }\frac{\frac{\varphi ^{n+2}}{(\varphi -\chi )(\varphi -\psi )}+\frac{%
\chi ^{n+2}}{(\chi -\varphi )(\chi -\psi )}+\frac{\psi ^{n+2}}{(\psi
-\varphi )(\psi -\chi )}}{\frac{\varphi ^{n+1}}{(\varphi -\chi )(\varphi
-\psi )}+\frac{\chi ^{n+1}}{(\chi -\varphi )(\chi -\psi )}+\frac{\psi ^{n+1}%
}{(\psi -\varphi )(\psi -\chi )}} \\
&=&\lim_{n\rightarrow \infty }\frac{\varphi ^{n+1}\left( \frac{\varphi }{%
(\varphi -\chi )(\varphi -\psi )}+\frac{\chi }{(\chi -\varphi )(\chi -\psi )}%
\frac{\chi ^{n+1}}{\varphi ^{n+1}}+\frac{\psi }{(\psi -\varphi )(\psi -\chi )%
}\frac{\psi ^{n+1}}{\varphi ^{n+1}}\right) }{\varphi ^{n}\left( \frac{%
\varphi }{(\varphi -\chi )(\varphi -\psi )}+\frac{\chi }{(\chi -\varphi
)(\chi -\psi )}\frac{\chi ^{n}}{\varphi ^{n}}+\frac{\psi }{(\psi -\varphi
)(\psi -\chi )}\frac{\psi ^{n}}{\varphi ^{n}}\right) } \\
&=&\lim_{n\rightarrow \infty }\frac{\varphi ^{n+1}}{\varphi ^{n}}\frac{%
\left( \frac{\varphi }{(\varphi -\chi )(\varphi -\psi )}+\frac{\chi }{(\chi
-\varphi )(\chi -\psi )}\left( \frac{\chi }{\varphi }\right) ^{n+1}+\frac{%
\psi }{(\psi -\varphi )(\psi -\chi )}\left( \frac{\psi }{\varphi }\right)
^{n+1}\right) }{\left( \frac{\varphi }{(\varphi -\chi )(\varphi -\psi )}+%
\frac{\chi }{(\chi -\varphi )(\chi -\psi )}\left( \frac{\chi }{\varphi }%
\right) ^{n}+\frac{\psi }{(\psi -\varphi )(\psi -\chi )}\left( \frac{\psi }{%
\varphi }\right) ^{n}\right) } \\
&=&\varphi \text{.}
\end{eqnarray*}

\section{Main Results}

In this section, we present our main results related to the difference
equation (\ref{equ:gutreb}). Our aim is to investigate the general solution
in explicit form of Eq.(\ref{equ:gutreb}) and the asymptotic behavior of
solutions of Eq.(\ref{equ:gutreb}).

\begin{theorem}
\label{theorem:aaaqssz}Let $\left \{ x_{n}\right \} _{n=-1}^{\infty }$ be a
solution of Eq.(\ref{equ:gutreb}). Then, for $n=0,1,2,...$, the form of
solutions $\left \{ x_{n}\right \} _{n=-1}^{\infty }$ is given by%
\begin{equation}
x_{n}=x_{n}=\frac{tV_{n-1}x_{-1}x_{0}+\left( V_{n+1}-rV_{n}\right)
x_{0}+V_{n}}{tV_{n}x_{-1}x_{0}+\left( V_{n+2}-rV_{n+1}\right) x_{0}+V_{n+1}},
\label{equ:mmmmmmn}
\end{equation}%
where $V_{n}$ is the $n$th generalized-Tribonacci number and the initial
conditions $x_{-1}$, $x_{0}\in 
\mathbb{R}
-F$, with $F$ is the forbidden set of Eq.(\ref{equ:gutreb}) given by%
\begin{equation*}
F=\dbigcup \limits_{n=-1}^{\infty }\left \{ \left( x_{-1},x_{0}\right)
:tV_{n}x_{-1}x_{0}+\left( V_{n+2}-rV_{n+1}\right) x_{0}+V_{n+1}=0\right \} 
\text{.}
\end{equation*}
\end{theorem}

\textbf{Proof.} First, by using the change of variables%
\begin{equation}
x_{n}=\frac{w_{n-1}}{w_{n}}\text{,}  \label{equ:ilkjhg}
\end{equation}%
Eq.(\ref{equ:gutreb}) is reduced to linear third order difference equation%
\begin{equation*}
w_{n+1}=\frac{\beta }{\gamma }w_{n}+\frac{\alpha }{\gamma }w_{n-1}+\frac{1}{%
\gamma }w_{n-2}\text{.}
\end{equation*}%
Set%
\begin{equation*}
r:=\frac{\beta }{\gamma }\text{, }s:=\frac{\alpha }{\gamma }\text{, \ }t:=%
\frac{1}{\gamma }\text{,}
\end{equation*}%
so we have%
\begin{equation*}
w_{n+1}=rw_{n}+sw_{n-1}+tw_{n-2}\text{.}
\end{equation*}%
Now, as done in [\ref{Stevic2018}], we describe initial values of three
sequences which will be repetitively defined and used in the rest of the
proof. Let%
\begin{equation*}
a_{1}:=r\text{, \ }b_{1}:=s\text{, \ }c_{1}:=t\text{.}
\end{equation*}%
We use an recurrent (iterative) method. Thus, we get%
\begin{eqnarray}
w_{n} &=&a_{1}w_{n-1}+b_{1}w_{n-2}+c_{1}w_{n-3}  \notag \\
&=&a_{1}\left( rw_{n-2}+sw_{n-3}+tw_{n-4}\right) +b_{1}w_{n-2}+c_{1}w_{n-3} 
\notag \\
&=&\left( ra_{1}+b_{1}\right) w_{n-2}+\left( sa_{1}+c_{1}\right)
w_{n-3}+ta_{1}w_{n-4}  \notag \\
&=&a_{2}w_{n-2}+b_{2}w_{n-3}+c_{2}w_{n-4}\text{,}  \label{equ:olkiu}
\end{eqnarray}%
where%
\begin{equation}
a_{2}:=ra_{1}+b_{1}\text{, \ }b_{2}:=sa_{1}+c_{1}\text{, \ }c_{2}:=ta_{1}%
\text{.}  \label{equ:polo}
\end{equation}%
By continuing iteration, it implies that%
\begin{eqnarray}
w_{n} &=&a_{2}w_{n-2}+b_{2}w_{n-3}+c_{2}w_{n-4}  \notag \\
&=&a_{2}\left( rw_{n-3}+sw_{n-4}+tw_{n-5}\right) +b_{2}w_{n-3}+c_{2}w_{n-4} 
\notag \\
&=&\left( ra_{2}+b_{2}\right) w_{n-3}+\left( sa_{2}+c_{2}\right)
w_{n-4}+ta_{2}w_{n-5}  \notag \\
&=&a_{3}w_{n-3}+b_{3}w_{n-4}+c_{3}w_{n-5}\text{,}  \label{eqt:sbsbtr}
\end{eqnarray}%
where%
\begin{equation}
a_{3}:=ra_{2}+b_{2}\text{, \ }b_{3}:=sa_{2}+c_{2}\text{, \ }c_{3}:=ta_{2}%
\text{.}  \label{equ:yhbngqasz}
\end{equation}%
Based on relations (\ref{equ:olkiu})-(\ref{equ:yhbngqasz}), we suppose that
for some $k\in 
\mathbb{N}
$ such that $2\leq k\leq n-1$, we have%
\begin{equation}
w_{n}=a_{k}w_{n-k}+b_{k}w_{n-k-1}+c_{k}w_{n-k-2}\text{,}
\label{equ:yhutrqxac}
\end{equation}%
and%
\begin{equation}
a_{k}:=ra_{k-1}+b_{k-1}\text{, \ }b_{k}:=sa_{k-1}+c_{k-1}\text{, \ }%
c_{k}:=ta_{k-1}\text{.}  \label{equ:rewqas}
\end{equation}%
Next, by continuing iteration, it follows that%
\begin{eqnarray*}
w_{n} &=&a_{k}w_{n-k}+b_{k}w_{n-k-1}+c_{k}w_{n-k-2} \\
&=&a_{k}\left( rw_{n-k-1}+sw_{n-k-2}+tw_{n-k-3}\right)
+b_{k}w_{n-k-1}+c_{k}w_{n-k-2} \\
&=&\left( ra_{k}+b_{k}\right) w_{n-k-1}+\left( sa_{k}+c_{k}\right)
w_{n-k-2}+ta_{k}w_{n-k-3} \\
&=&a_{k+1}w_{n-k-1}+b_{k+1}w_{n-k-2}+c_{k+1}w_{n-k-3}\text{,}
\end{eqnarray*}%
where%
\begin{equation*}
a_{k+1}:=ra_{k}+b_{k}\text{, \ }b_{k+1}:=sa_{k}+c_{k}\text{, \ }%
c_{k+1}:=ta_{k}\text{.}
\end{equation*}

Now, we maintain sequences $a_{k}$, $b_{k}$ and $c_{k}$ for some nonpositive
values of index $k$. Notice that since $\gamma \neq 0$, the recurrent
relations in (\ref{equ:rewqas}) can be really used for computing values of
sequences $a_{k}$, $b_{k}$ and $c_{k}$ for every $k\leq 0$.

Using the recurrent relations with the indices $k=1$, $k=0$ and $k=-1$,
respectively, after some computations, it implies that%
\begin{eqnarray*}
a_{0} &=&\frac{c_{1}}{c}=1 \\
b_{0} &=&a_{1}-aa_{0}=a-a.1=0 \\
c_{0} &=&b_{1}-ba_{0}=b-b.1=0 \\
a_{-1} &=&\frac{c_{0}}{c}=0 \\
b_{-1} &=&a_{0}-aa_{-1}=1-a.0=1 \\
c_{-1} &=&b_{0}-ba_{-1}=0-b.0=0 \\
a_{-2} &=&\frac{c_{-1}}{c}=0 \\
b_{-2} &=&a_{-1}-aa_{-2}=0-a.0=0 \\
c_{-2} &=&b_{-1}-ba_{-2}=1-b.0=1\text{.}
\end{eqnarray*}%
Thus, we obtain%
\begin{equation}
\begin{array}{ccc}
a_{0}=1 & a_{-1}=0 & a_{-2}=0 \\ 
b_{0}=0 & b_{-1}=1 & b_{-2}=0 \\ 
c_{0}=0 & c_{-1}=0 & c_{-2}=1\text{.}%
\end{array}
\label{iv:bbbbiiii}
\end{equation}%
From (\ref{equ:rewqas}), we get%
\begin{equation}
a_{n}=ra_{n-1}+sa_{n-2}+ta_{n-3}\text{,}  \label{equ:wqsxiopk}
\end{equation}%
\begin{equation}
b_{n}=a_{n+1}-ra_{n}\text{,}  \label{equ:tfvcds}
\end{equation}%
\begin{equation}
c_{n}=ta_{n-1}\text{,}  \label{equ:okjhfedsq}
\end{equation}%
for $n\in 
\mathbb{N}
$.

If we get $k=n$ in (\ref{equ:yhutrqxac}), we have%
\begin{equation*}
w_{n}=a_{n}w_{0}+b_{n}w_{-1}+c_{n}w_{-2}\text{,}
\end{equation*}%
for $n\in 
\mathbb{N}
_{0}$.

From (\ref{equ:wqsxiopk})-(\ref{equ:okjhfedsq}), we obtain%
\begin{equation}
w_{n}=a_{n}w_{0}+\left( a_{n+1}-ra_{n}\right) w_{-1}+ta_{n-1}w_{-2}\text{,}
\label{equ:ygrew}
\end{equation}%
for $n\in 
\mathbb{N}
_{0}$

Using (\ref{equ:ygrew}) in (\ref{equ:ilkjhg}), we get%
\begin{equation*}
x_{n}=\frac{a_{n-1}w_{0}+\left( a_{n}-ra_{n-1}\right) w_{-1}+ta_{n-2}w_{-2}}{%
a_{n}w_{0}+\left( a_{n+1}-ra_{n}\right) w_{-1}+ta_{n-1}w_{-2}}\text{,}
\end{equation*}%
it follows that%
\begin{equation*}
x_{n}=\frac{ta_{n-2}x_{-1}x_{0}+\left( a_{n}-ra_{n-1}\right) x_{0}+a_{n-1}}{%
ta_{n-1}x_{-1}x_{0}+\left( a_{n+1}-ra_{n}\right) x_{0}+a_{n}}
\end{equation*}%
or equivalently%
\begin{equation*}
x_{n}=\frac{ta_{n-2}x_{-1}x_{0}+\left( a_{n}-ra_{n-1}\right) x_{0}+a_{n-1}}{%
ta_{n-1}x_{-1}x_{0}+\left( sa_{n-1}+ta_{n-2}\right) x_{0}+a_{n}}\text{.}
\end{equation*}%
From initial values (\ref{iv:bbbbiiii}) and definitions of sequences $a_{n}$
and $V_{n}$, we have%
\begin{equation*}
a_{n}=V_{n+1}\text{,}
\end{equation*}%
with the backward shifted initial values of the sequence $a_{n}$. Then, it
follows%
\begin{equation*}
x_{n}=\frac{tV_{n-1}x_{-1}x_{0}+\left( V_{n+1}-rV_{n}\right) x_{0}+V_{n}}{%
tV_{n}x_{-1}x_{0}+\left( V_{n+2}-rV_{n+1}\right) x_{0}+V_{n+1}}\text{,}
\end{equation*}%
or%
\begin{equation*}
x_{n}=\frac{tV_{n-1}x_{-1}x_{0}+\left( V_{n+1}-rV_{n}\right) x_{0}+V_{n}}{%
tV_{n}x_{-1}x_{0}+\left( sV_{n}+tV_{n-1}\right) x_{0}+V_{n+1}}\text{.}
\end{equation*}%
The proof is complete.

Now, we will analyze five special cases of the above theorem according to
the states of $r$, $s$, $t$.

\textbf{Case 1: }$r=s=t=1$

In this case, the $\left( a_{n}\right) $ sequence has the following
recurrence relation%
\begin{equation*}
a_{n}=a_{n-1}+a_{n-2}+a_{n-3}\text{,}
\end{equation*}%
such that a few terms of this sequence are%
\begin{equation}
a_{-2}=0,a_{-1}=0,a_{0}=1,a_{1}=1,a_{2}=2,a_{3}=4\text{.}
\label{iv:eqwsamlip}
\end{equation}

Then, from initial values (\ref{iv:eqwsamlip}) and definitions of sequences $%
a_{n}$ and $T_{n}$ which is Tribonacci numbers, we have%
\begin{equation*}
a_{n}=T_{n+1}\text{,}
\end{equation*}%
with the backward shifted initial values of the sequence $a_{n}$.

Hence, we obtain%
\begin{equation*}
x_{n}=\frac{T_{n-1}x_{-1}x_{0}+\left( T_{n+1}-T_{n}\right) x_{0}+T_{n}}{%
T_{n}x_{-1}x_{0}+\left( T_{n}+T_{n-1}\right) x_{0}+T_{n+1}}.
\end{equation*}%
\textbf{Case 2: }$r=0$, $s=t=1$

In this case, the $\left( a_{n}\right) $ sequence has the following
recurrence relation%
\begin{equation*}
a_{n}=a_{n-2}+a_{n-3}\text{,}
\end{equation*}%
such that a few terms of this sequence are%
\begin{equation}
a_{-2}=0,\text{ }a_{-1}=0,\text{ }a_{0}=1,\text{ }a_{1}=0,\text{ }a_{2}=1,%
\text{ }a_{3}=1\text{, }a_{4}=1\text{.}  \label{iv:bbbbbbbbs}
\end{equation}

Then, from initial values (\ref{iv:bbbbbbbbs}) and definitions of sequences $%
a_{n}$ and $P_{n}$ which is Padovan numbers, we get%
\begin{equation*}
a_{n+2}=P_{n}\text{,}
\end{equation*}%
with the forward shifted initial values of the sequence $a_{n}$.

Therefore, we have%
\begin{equation*}
x_{n}=\frac{P_{n-4}x_{-1}x_{0}+P_{n-2}x_{0}+P_{n-3}}{%
P_{n-3}x_{-1}x_{0}+P_{n-1}x_{0}+P_{n-2}}\text{.}
\end{equation*}%
\textbf{Case 3: }$r=0$, $s=t=1$

In this case, the $\left( a_{n}\right) $ sequence has the following
recurrence relation%
\begin{equation*}
a_{n}=a_{n-2}+a_{n-3}\text{,}
\end{equation*}%
such that a few terms of this sequence are%
\begin{equation}
a_{-2}=0,\text{ }a_{-1}=0,\text{ }a_{0}=1,\text{ }a_{1}=0,\text{ }a_{2}=1,%
\text{ }a_{3}=1\text{, }a_{4}=1\text{.}  \label{iv:ujnmb}
\end{equation}

Then, from initial values (\ref{iv:ujnmb}) and definitions of sequences $%
a_{n}$ and $S_{n}$ which is Padovan-Perrin numbers, we have%
\begin{equation*}
a_{n}=S_{n+2}\text{,}
\end{equation*}%
with the backward shifted initial values of the sequence $a_{n}$.

Thus, we obtain%
\begin{equation*}
x_{n}=\frac{S_{n}x_{-1}x_{0}+S_{n+2}x_{0}+S_{n+1}}{%
S_{n+1}x_{-1}x_{0}+S_{n+3}x_{0}+S_{n+2}}\text{.}
\end{equation*}%
\textbf{Case 4: }$r=1$, $s=0$, $t=1$

In this case, the $\left( a_{n}\right) $ sequence has the following
recurrence relation%
\begin{equation*}
a_{n}=a_{n-1}+a_{n-3}\text{,}
\end{equation*}%
such that a few terms of this sequence are%
\begin{equation}
a_{-2}=0,\text{ }a_{-1}=0,\text{ }a_{0}=1,\text{ }a_{1}=1,\text{ }a_{2}=1,%
\text{ }a_{3}=2\text{, }a_{4}=3\text{.}  \label{iv:iiiiiiiiiiiipo}
\end{equation}

Then, from initial values (\ref{iv:iiiiiiiiiiiipo}) and definitions of
sequences $a_{n}$ and $N_{n}$ which is Narayana numbers, we have%
\begin{equation*}
a_{n}=N_{n+1}\text{,}
\end{equation*}%
with the backward shifted initial values of the sequence $a_{n}$.

Fromhere, we have%
\begin{equation*}
x_{n}=\frac{N_{n-1}x_{-1}x_{0}+N_{n-2}x_{0}+N_{n}}{%
N_{n}x_{-1}x_{0}+N_{n-1}x_{0}+N_{n+1}}\text{.}
\end{equation*}%
\textbf{Case 5: }$r=s=1$, $t=2$

In this case, the $\left( a_{n}\right) $ sequence has the following
recurrence relation%
\begin{equation*}
a_{n}=a_{n-1}+a_{n-2}+2a_{n-3}\text{,}
\end{equation*}%
such that a few terms of this sequence are%
\begin{equation}
a_{-2}=0,\text{ }a_{-1}=0,\text{ }a_{0}=1,\text{ }a_{1}=1,\text{ }a_{2}=2,%
\text{ }a_{3}=5\text{, }a_{4}=9\text{.}  \label{iv:aemnb}
\end{equation}

Next, from initial values (\ref{iv:aemnb}) and definitions of sequences $%
a_{n}$ and $J_{n}$ which is third order Jacobsthal numbers, we have%
\begin{equation*}
a_{n}=J_{n+1}\text{,}
\end{equation*}%
with the backward shifted initial values of the sequence $a_{n}$.

Herefrom, we get%
\begin{equation*}
x_{n}=\frac{2J_{n-1}x_{-1}x_{0}+\left( J_{n+1}-J_{n}\right) x_{0}+J_{n}}{%
2J_{n}x_{-1}x_{0}+\left( J_{n+2}-J_{n+1}\right) x_{0}+J_{n+1}}\text{.}
\end{equation*}

\begin{theorem}
\label{theorem:brtgyz}Eq.(\ref{equ:gutreb}) has unique equilibrium point $%
\overline{x}=\theta $ and $\theta $ is locally asymptotically stable.
\end{theorem}

\textbf{Proof.} Equilibrium point of Eq.(\ref{equ:gutreb}) is the real roots
of the equation%
\begin{equation*}
\overline{x}=\frac{\gamma }{\overline{x}\left( \overline{x}+\alpha \right)
+\beta }\text{.}
\end{equation*}%
After simplification, we get the following cubic equation%
\begin{equation}
\overline{x}^{3}+\alpha \overline{x}^{2}+\beta \overline{x}-\gamma =0\text{.}
\label{equ:zzzaaassa}
\end{equation}%
Then, the roots of the cubic equation (\ref{equ:zzzaaassa}) are given by%
\begin{eqnarray*}
\mu &=&\mu \left( \alpha ,\beta ,\gamma \right) =-\frac{\alpha }{3}+C+D\text{%
,} \\
\sigma &=&\sigma \left( \alpha ,\beta ,\gamma \right) =-\frac{\alpha }{3}%
+\omega C+\omega ^{2}D\text{,} \\
\phi &=&\phi \left( \alpha ,\beta ,\gamma \right) =-\frac{\alpha }{3}+\omega
^{2}C+\omega D\text{,}
\end{eqnarray*}%
where%
\begin{eqnarray*}
C &=&\left( \frac{-\alpha ^{3}}{27}+\frac{\alpha \beta }{6}+\frac{\gamma }{2}%
+\sqrt{\Delta }\right) ^{1/3}\text{, }D=\left( \frac{-\alpha ^{3}}{27}+\frac{%
\alpha \beta }{6}+\frac{\gamma }{2}-\sqrt{\Delta }\right) ^{1/3} \\
\Delta &=&\Delta (r,s,t)=-\frac{\alpha ^{3}\gamma }{27}-\frac{\alpha
^{2}\beta ^{2}}{108}+\frac{\alpha \beta \gamma }{6}+\frac{\beta ^{3}}{27}+%
\frac{\gamma ^{2}}{4}
\end{eqnarray*}%
and%
\begin{equation*}
\omega =\frac{-1+i\sqrt{3}}{2}=\exp \left( 2\pi i/3\right)
\end{equation*}%
is a primitive cube root of unity.\ So, the root $\mu $ is only real number.
So, the unique equilibrium point of Eq.(\ref{equ:gutreb}) is $\overline{x}%
=\mu $.

Now, we demonstrate that the equilibrium point of Eq.(\ref{equ:gutreb}) is
locally asymptotically stable.

Let $I$ be an interval of real numbers and consider the function%
\begin{equation*}
f:I^{2}\rightarrow I
\end{equation*}%
defined by%
\begin{equation*}
f\left( x,y\right) =\frac{\gamma }{x\left( y+\alpha \right) +\beta }\text{.}
\end{equation*}%
The linearized equation of Eq.(\ref{equ:gutreb}) about the equilibrium point 
$\overline{x}=\mu $ is%
\begin{equation*}
z_{n+1}=pz_{n}+qz_{n-1}\text{,}
\end{equation*}%
where%
\begin{eqnarray*}
p &=&\frac{\partial f\left( \overline{x},\overline{x}\right) }{\partial x}=%
\frac{-\gamma \left( \mu +\alpha \right) }{\left( \mu \left( \mu +\alpha
\right) +\beta \right) ^{2}} \\
&=&\frac{-\gamma \left( \mu +\alpha \right) }{\left( \mu ^{2}+\alpha \mu
+\beta \right) ^{2}} \\
&=&\frac{-\gamma \left( \mu +\alpha \right) }{\left( \frac{\gamma }{\mu }%
\right) ^{2}} \\
&=&\frac{-\left( \mu ^{3}+\alpha \mu ^{2}\right) }{\gamma } \\
&=&\frac{\beta \mu -\gamma }{\gamma }\text{,} \\
q &=&\frac{\partial f\left( \overline{x},\overline{x}\right) }{\partial y}=%
\frac{-\gamma \mu }{\left( \mu \left( \mu +\alpha \right) +\beta \right) ^{2}%
} \\
&=&\frac{-\gamma \mu }{\left( \mu ^{2}+\alpha \mu +\beta \right) ^{2}} \\
&=&\frac{-\gamma \mu }{\left( \frac{\gamma }{\mu }\right) ^{2}} \\
&=&-\frac{\mu ^{3}}{\gamma }\text{,}
\end{eqnarray*}%
and the corresponding characteristic equation is%
\begin{equation*}
\lambda ^{2}-\left( \frac{\beta \mu -\gamma }{\gamma }\right) \lambda +\frac{%
\mu ^{3}}{\gamma }=0\text{.}
\end{equation*}%
Consider two functions defined by%
\begin{equation*}
a\left( \lambda \right) =\lambda ^{2}\text{, \ }b\left( \lambda \right)
=\left( \frac{\beta \mu -\gamma }{\gamma }\right) \lambda -\frac{\mu ^{3}}{%
\gamma }\text{.}
\end{equation*}%
We have%
\begin{equation*}
\left\vert \frac{\beta \mu -\gamma }{\gamma }-\frac{\mu ^{3}}{\gamma }%
\right\vert <1\text{.}
\end{equation*}%
Then, it follows that%
\begin{equation*}
\left\vert b\left( \lambda \right) \right\vert <\left\vert a\left( \lambda
\right) \right\vert \text{, \ for all }\lambda :\left\vert \lambda
\right\vert =1\text{.}
\end{equation*}%
Therefore, by Rouche's Theorem, all zeros of $P\left( \lambda \right)
=a\left( \lambda \right) -b\left( \lambda \right) =0$ lie in $\left\vert
\lambda \right\vert <1$. Hereby, by Theorem (\ref{theorem:lkjhgfvcnm}), we
have that the unique equilibrium point of Eq.(\ref{equ:gutreb}) $\overline{x}%
=\mu $ is locally asymptotically stable.

\begin{theorem}
\label{theo:iiiiiiiiiiiiiiiiiiiibbbb}Assume that $\mu \varphi =1$. Then, the
equilibrium point of Eq.(\ref{equ:gutreb}) is globally asymptotically stable.
\end{theorem}

\textbf{Proof.} Let $\left\{ x_{n}\right\} _{n\geq -1}$ be a solution of Eq.(%
\ref{equ:gutreb}). By Theorem (\ref{theorem:brtgyz}), we need only to prove
that the equilibrium point $\mu $ is global attractor, that is%
\begin{equation*}
\lim_{n\rightarrow \infty }x_{n}=\mu \text{.}
\end{equation*}%
From Theorem (\ref{theorem:aaaqssz}), (\ref{chaequ:jhndf}) and (\ref%
{lim:ygbhf}), it follows that%
\begin{eqnarray*}
\lim_{n\rightarrow \infty }x_{n} &=&\lim_{n\rightarrow \infty }\frac{%
tV_{n-1}x_{-1}x_{0}+\left( V_{n+1}-rV_{n}\right) x_{0}+V_{n}}{%
tV_{n}x_{-1}x_{0}+\left( V_{n+2}-rV_{n+1}\right) x_{0}+V_{n+1}} \\
&=&\lim_{n\rightarrow \infty }\frac{tV_{n-1}\left( x_{-1}x_{0}+\left( \frac{1%
}{t}\frac{V_{n+1}}{V_{n-1}}-\frac{r}{t}\frac{V_{n}}{V_{n-1}}\right) x_{0}+%
\frac{V_{n}}{V_{n-1}}\right) }{tV_{n}\left( x_{-1}x_{0}+\left( \frac{1}{t}%
\frac{V_{n+2}}{V_{n}}-\frac{r}{t}\frac{V_{n+1}}{V_{n}}\right) x_{0}+\frac{%
V_{n+1}}{V_{n}}\right) } \\
&=&\left( \frac{x_{-1}x_{0}+\left( \frac{1}{t}\varphi ^{2}-\frac{r}{t}%
\varphi \right) x_{0}+\varphi }{x_{-1}x_{0}+\left( \frac{1}{t}\varphi ^{2}-%
\frac{r}{t}\varphi \right) x_{0}+\varphi }\right) \lim_{n\rightarrow \infty }%
\frac{V_{n-1}}{V_{n}} \\
&=&\lim_{n\rightarrow \infty }\frac{V_{n-1}}{V_{n}} \\
&=&\frac{1}{\varphi } \\
&=&\mu \text{.}
\end{eqnarray*}%
This completes the proof.

Note that when $\alpha =\beta =\gamma =1$, our assumption in Theorem \ref%
{theo:iiiiiiiiiiiiiiiiiiiibbbb} is immediately seen. Indeed,

$\mu \varphi =\allowbreak \frac{1}{3\gamma }\left( \frac{1}{3\sqrt[3]{2}}%
\sqrt[3]{27\gamma -2\alpha ^{3}+9\alpha \beta -3\sqrt{3}S}-\frac{1}{3}\alpha
+\frac{1}{3\sqrt[3]{2}}\sqrt[3]{27\gamma -2\alpha ^{3}+9\alpha \beta +3\sqrt{%
3}S}\right) $

$\allowbreak \left( \beta +\frac{1}{\sqrt[3]{2}}\sqrt[3]{2\beta
^{3}+27\gamma ^{2}-3\sqrt{3}\gamma S+9\alpha \beta \gamma }+\frac{1}{\sqrt[3]%
{2}}\sqrt[3]{2\beta ^{3}+27\gamma ^{2}+3\sqrt{3}\gamma S+9\alpha \beta
\gamma }\right) $

$S=\sqrt{-4\alpha ^{3}\gamma -\alpha ^{2}\beta ^{2}+18\alpha \beta \gamma
+4\beta ^{3}+27\gamma ^{2}}$.

Then, in the case $\alpha =\beta =\gamma =1$, it follows that $\mu \varphi
=1 $.

\end{document}